\documentclass[a4paper,draft, 10pt]{amsproc}
\usepackage{amssymb}
\usepackage{amscd} 
\usepackage{enumerate, multicol}
\usepackage{cite}
\usepackage[autostyle]{csquotes}

\theoremstyle{plain}
 \newtheorem{thm}{Theorem}[section]

\theoremstyle{definition}

\theoremstyle{remark}
 
 \numberwithin{equation}{section}

\renewcommand{\leq}{\leqslant}

\title[Darboux's Theorem]{Another Proof of Darboux's Theorem}



\author[M.Bhandari]{\bfseries Dr. Mukta Bhandari}

\address{
Department of Mathematics \\ 
Chowan University   \\ 
Murfreesboro\\
USA}
\email{bhandm@chowan.edu}

\begin{document}

\vspace{18mm} \setcounter{page}{1} \thispagestyle{empty}

\begin{abstract}
We know that a continuous function on a closed interval satisfies the Intermediate Value Property. Likewise, the derivative function of a differentiable function on a closed interval satisfies the IVP property which is known as the Darboux's theorem in any real analysis course. Most of the proofs found in the literature use the Extreme Value Property of a continuous function. In this paper, I am going to present a simple and elegant proof of the Darboux's theorem using the Intermediate Value Theorem and the Rolles theorem
\end{abstract}

\maketitle

\section{Some Preliminary Background and Known Proofs}  

In this section we state the Darboux's theorem and give the known proofs from various literatures.

A function $f\colon[a,b]\rightarrow \mathbb{R}$ is said to satify the intermediate value property on $[a,b]$ if for every $\lambda $ between $f(a)$ and $f(b)$, there exists $c\in (a,b)$ such that $f(c) = \lambda$.
We know, from intermediate value theorem, that a function  $f\colon [a,b]\rightarrow \mathbb{R}$ continuous on $[a,b]$ satisfies the intermediate value property on $[a,b]$. That is

$$ \text{continuous function}\quad \Rightarrow\quad \text{Intermediate Value Property}.$$

However, the converse of the intermediate value theorem is not necessarily true. The function $f\colon \mathbb{R} \rightarrow \mathbb{R}$ defined by

\begin{displaymath}
   f(t) = \left\{
     \begin{array}{lr}
       \sin\left( \dfrac{1}{t}\right)  & for \quad t\neq 0,\\
       0 & for \quad t=0, \\
     \end{array}
   \right.
\end{displaymath}
which is discontinuous at $t=0$, provides a counter example. Indeed, the French mathematician Darboux, proves in $1875$ that derivative function of a differentiable function satisfies the intermediate value property ~\cite{darboux}. Darboux also provides an example of a differentiable function with discontinuous derivative.

The statement of the Darboux's theorem follows here.

\begin{thm}[Darboux's Theorem]\label{thm:darbouxtheorem}
If $f$ is differentiable on $[a,b]$ and if $\lambda$ is a number between $f^{\prime}(a)$ and $f^{\prime}(b)$, then there is at least one point $c\in (a,b)$ such that $f^{\prime}(c) = \lambda$.
\end{thm}

 This indicates clearly that the converse of the intermediate value theorem is not necessarily true because derivative function of a differentiable function is not necessarily continuous. For example

\begin{displaymath}
   f(t) = \left\{
     \begin{array}{lr}
       t^2\sin\left( \dfrac{1}{t}\right)  & for \quad t\neq 0,\\
       0 & for \quad t=0 \\
     \end{array}
   \right.
\end{displaymath}
is continuous and has discontinuous derivative which satisfies the intermediate value property by the Darboux's theorem \ref{thm:darbouxtheorem}. Thus Draboux's provides a huge class of functions satisfying the intermediate value property, namely the derivative functions of differentiable functions.
The main purpose of this short note is to provide an alternate proof of this theorem different from the standard proofs found in many textbooks on real analysis.

The standard proofs found in various literatures goes as follows:

\begin{proof}
Suppose that $f^{\prime}(a) < \lambda < f^{\prime}(b)$. Let $F:[a,b]\rightarrow \mathbb{R}$ be defined by $F(x) = f(x) -\lambda x$ so that $F^{\prime}(x) = f^{\prime}(x) -\lambda $. Then $F$ is differentiable on $[a,b]$ because so is the function $f$ by hypothesis. We find  $F^{\prime}(a) = f^{\prime}(a) -\lambda  <0$ and  $F^{\prime}(b) = f^{\prime}(b) -\lambda >0$. Note that $F^{\prime}(a) <0$ means $F(t_1)< F(a)$ for some $t_1\in (a,b)$. Also, for $F^{\prime}(b)>0$, we can find $t_2\in (a,b)$ such that $F(t_2) < F(b)$. Thus neither $a$ nor $b$ can be a point where $F$ attains absolute minimum. Since $F$ is continuous on $[a,b]$, it must attain its relative minimum at some point $c\in (a,b)$ by Extreme value theorem for a continuous function. This means that $F^{\prime}(c) = 0$  by Fermat's theorem  and therefore $f^{\prime}(c) = \lambda$ as desired. The proof follows by similar arguments if  $f^{\prime}(b) < \lambda < f^{\prime}(a)$.
\end{proof}

The above proof can be found in various textbooks of undergraduate level real analysis course including W. Rudin ~\cite{walterrudin}, M. Spivak ~\cite{mspivak}, Bartle and Sherbert ~\cite{bartlesherbert}, K. Ross ~\cite{kaross}, W. R. Wade ~\cite{wade}, J. M. Howe ~\cite{howie}, T. M. Apostol ~\cite{apostol}, R. Boas ~\cite{boas}, G. Darboux ~\cite{darboux}, S. Krantz ~\cite{krantz}. It is left as an exercise with some hints in some of the textbooks while others give incomplete proofs leaving readers to fill in the gaps. The main trouble is in the argument that the function $F$ can attain minimum neither at $x=a$ nor at $x=b$. Some authors provide the $\epsilon-\delta$ proof for this part of the argument. Lars Oslen, in his paper ~\cite{oslen}, makes the following comment:

\enquote{Most students typically either think that this is obvious (and that the lecture is being over pedantic by insisting on a proof), or they see the need for a proof but find the $\epsilon-\delta$-gymnastics in the proof less than convincing.}

Lars Oslen then provides an new proof of the Darboux's theorem based only on the Mean Value theorem  for a differentiable function and intermediate value theorem  for a continuous function. The author, Lars Oslen, claims that his proof is more convincing than the standard proofs found in many textbooks. I am going to present his proof as it is from his paper ~\cite{oslen}.

\medskip

\noindent{\bf Proof of Darboux's Theorem by Lars Oslen ~\cite{oslen}:}

\medskip
We may clearly assume that $y$ lies stricly between $f^{\prime}(a)$ and $f^{\prime}(b)$. Define continuous functions $f_a, f_b\colon [a,b]\rightarrow \mathbb{R}$ by

\begin{displaymath}
   f_a(t) = \left\{
     \begin{array}{lr}
       f^{\prime}(a) & for \quad t=a,\\
       \dfrac{f(a)-f(t)}{a-t} & for \quad t\neq a, \\
     \end{array}
   \right.
\end{displaymath}
and

\begin{displaymath}
   f_b(t) = \left\{
     \begin{array}{lr}
       f^{\prime}(b) & for \quad t=b,\\
       \dfrac{f(t)-f(b)}{t-b} & for \quad t\neq b. \\
     \end{array}
   \right.
\end{displaymath}

It follows that $f_a(a) = f^{\prime}(a)$, $f_a(b)=f_b(a)$, and $f_b(b) = f^{\prime}(b)$. Hence, $y$ lies between $f_a(a)$ and $f_a(b)$, or $y$ lies between $f_b(a)$ and $f_b(b)$.
\par
If $y$ lies between $f_a(a)$ and $f_a(b)$, then (by continuity of $f_a$) there exists $s$ in $(a,b]$ with

$$y= f_a(s) = \dfrac{f(s) - f(a)}{s-a}.$$
The mean value theorem ensures that there exists $x$ in $[a,s]$ such that

$$y = \dfrac{f(s)-f(a)}{s-a}=f^{\prime}(x).$$

If $y$ lies between $f_b(a)$ and $f_b(b)$, then an analogous argument (exploiting the continuity of $f_b$) shows that there exists $s$ in $[a,b)$ such that

$$y = \dfrac{f(b)-f(s)}{b-s}=f^{\prime}(x).$$
This completes the proof.

This proof is more convincing as the author has claimed in his paper in the sense that it avoids the $\epsilon-\delta$-arguments.  Author has used the Carath\'{e}odory's theorem for the existence of the auxiliary functions $f_a$ and $f_b$. See T. M. Apostol ~\cite{apostol}, Bartle and Sherbert ~\cite{bartlesherbert} for example. Use of this theorem for the existence of the functions $f_a$ and $f_b$ may still be rather too much for an undergraduate students beginning real analysis course.
The proof I am going to give in the next section uses monotonicity property of a differentiable function, and the  two standard and familiar theorems   intermediate value theorem  and the Rolle's theorem from the first year calculus course. This also uses no $\epsilon-\delta$-arguments.

\section{Another Proof of Darboux's Theorem}

\begin{proof}

We assume that $ f^{\prime}(a) < \lambda < f^{\prime}(b)$, and consider a function $F\colon[a,b]\rightarrow \mathbb{R}$ defined by $F(x) = f(x) - \lambda x$.
Then $F^{\prime}(x) = f^{\prime}(x) - \lambda$, and find that $F^{\prime}(a) = f^{\prime}(a) - \lambda <0$ and $F^{\prime}(b)= f^{\prime}(b)-\lambda >0$.
This means that the function $F$ is not monotonic on $[a,b]$ in the sense that there exist $x, y, z\in [a,b]$ such that $x < y < z$ satisfying  exactly one and only one of the following conditions $(a)$ or $(b)$.

\medskip

\begin{enumerate}[(a)]
\item $F(x) < F(y)$ and $F(y) > F(z)$ We prove the theorem for this case.
\item $F(x) > F(y)$ and $F(y) < F(z)$. Proof in this case follows by a similar argument.
\end{enumerate}

\medskip

Suppose $(a)$ holds. Then we encounter three cases:

\begin{enumerate}[(i)]
\item  $F(x) < F(z)$
\item $F(z) < F(x)$
\item $F(z) = F(x)$.
\end{enumerate}
\par
\medskip

\noindent{Case(i):} Assume that $F(x) < F(z)$. Then $F(x) < F(z) < F(y)$. F is differentiable on $[a,b]$ because so is $f$ by hypothesis. That is $F$ is continuous as well on $[a,b]$. So the Intermediate Value theorem applies to $F$ on $[x,y]\subseteq [a,b]$ and we obtain $d\in (x, y) \subseteq (a,b)$ such that $F(d) = F(z)$.  Note that

$$a \leq x < d < y < z \leq b.$$

Now we can apply the Rolles theorem  to $F$ on the closed interval $[d,z]$ and obtain $c\in (d,z)\subseteq (x,z)\subseteq (a,b)$ such that $F^{\prime}(c) = 0$  which leads to $f^{\prime}(c) = \lambda$ as desired.
\par
\medskip

\noindent{Case (ii):} Assume that $F(z) < F(x)$. Then $F(z) < F(x) < F(y)$. The theorem then follows by arguments similar to that of case (i).

\medskip

\noindent{Case (iii):} Assume that $F(z) = F(x)$. Then we can use the Rolle's theorem  for $F$ directly on $[x,z]$ and obtain a $c\in (x,z)\subseteq (a,b)$ such that $F^{\prime}(c) = 0$ leading to $f^{\prime}(c)=\lambda$ as desired.

\par

 We arrive at the same conclusion if $(b)$ holds by similar  argument. This completes the proof of the theorem.

\end{proof}

\bibliographystyle{amsplain}

\end{document}